\documentclass[aps,pre,twocolumn,showpacs,floatfix,preprintnumbers,amsmath,amssymb]{revtex4-1}
\usepackage[english]{babel}
\usepackage[utf8]{inputenc}
\usepackage[colorinlistoftodos, color=green!40, prependcaption]{todonotes}
\usepackage{float}
\usepackage{mathtools}
\usepackage{tikz}
\usepackage{physics}
\usepackage{xcolor}
\usepackage{graphicx}
\usepackage{caption}
\usepackage{subcaption}
\usepackage{mathrsfs}
\usepackage{amssymb}
\usepackage{amsthm}
\usepackage{multirow}
\usepackage[hidelinks]{hyperref} 
\hypersetup{
  colorlinks   = true, 
  urlcolor     = blue, 
  linkcolor    = blue, 
  citecolor   = red 
}

\begin{document}
\title{Machine learning the interaction network in coupled dynamical systems}
\author{Pawan R. Bhure and M. S. Santhanam}
\affiliation{Physics Department, Indian Institute of Science Education and Research, Pune 411008, India}
\date{\today}

\begin{abstract}
The study of interacting dynamical systems continues to attract research interest in various fields of science and engineering. In a collection of interacting particles, the interaction network contains information about how various components interact with one another. Inferring the information about the interaction network from the dynamics of agents is a problem of long-standing interest. In this work, we employ a self-supervised neural network model to achieve two outcomes: to recover the interaction network and to predict the dynamics of individual agents. Both these information are inferred solely from the observed trajectory data. This work presents an application of the Neural Relational Inference model to two dynamical systems: coupled particles mediated by Hooke's law interaction and coupled phase (Kuramoto) oscillators.
\end{abstract}

\maketitle

\section{Introduction}
Interacting dynamical systems are ubiquitous, and one of the primary goals in physics and other areas of sciences and engineering is to obtain a deeper insight into the working of such systems \cite{barrat2008}. They appear in diverse physical settings ranging from a collection of gas molecules, interactions among human agents, vehicular traffic, economic interactions among agents to dynamics of neurons and other physiological processes \cite{vespignani2012}. Generally, they consist of several agents or entities that interact with one another, i.e. one agent influences the state of the other agents depending on how strongly they are coupled to one another. Hence, such coupled dynamical systems can be thought of constituting a network in which each node represents a particle (or, some basic dynamical unit of interest) and the edges represent the strength of interactions between the nodes \cite{battiston2020, aleta2019}. In the last two decades or more, extensive research investigations have been carried out to understand the emergent properties of such systems such as synchronization, extreme events, amplitude death and resilience of such networks \cite{barrat2008,bar2019}.

In these theoretical and modelling approaches to interacting systems, it is implicitly assumed that the connections between constituent units are well known, i.e, they are provided upfront as a part of the problem statement. This is the forward modelling approach that is usually employed which entails understanding the dynamics given the network parameters such as the connections between the nodes. However, in many situations the inverse problems is of practical interest \cite{li2021percolation}. A major challenge is that we often do not have access to the precise physical interactions or the laws governing the dynamics. Often, we have access only to the measured time series of some of the dynamical variables. For instance, in the case of vehicular traffic, it is known to be an interacting system though the form of interactions among vehicles is unknown. For any study of vehicular traffic, only a partial and coarse-grained record of vehicular time series is available. Similarly, in the context of stock markets, precise nature of interactions among stocks are not known though it can be inferred from the evolution of stock data and market indices.  Thus, a general scenario in many contexts can be posed as follows -- in a network of interacting systems for which we have access only to a time series record of evolving dynamics, can we use this information to infer the interactions among the units \cite{Peixoto2019}. This is the inverse problem of interest in this work.

Since the last decade, this problem has been extensively addressed by various authors, see Ref. \cite{TimCas2014} for an overview of various approaches. In general, many techniques exist to uncover communities and modularity in networks \cite{For2010}. It is usual for network edge detection to be based on observed time series data. Early proposals included approaches based on observation of stable response dynamics -- phase differences and collective frequency of constituent dynamical units -- in the presence of constant driving \cite{Tim2007}. Over the years, a host of approaches had been employed  -- dynamics based techniques \cite{GuoFu2010,ShaTim2011,RubMarBia2014,ChiLaiLeu2015,HanSheWan2015}, statistical techniques \cite{NapSau2008,KraEdeCas2009,MaZhaLai2017,MaCheLai2018,ShiSheJin2021}, compressive sensing based techniques \cite{WanLaiGre2011,SheWanWen2014,BarZho2019}, and also inference directly from empirical data \cite{TopEro2023,GaoYan2022,CheZhaZha2022,WuHaoLiu2022}.

The emergence of machine learning had led to works that have attempted to solve graph inference problem using tools drawn from neural networks \cite{CheLai2018,XiaSuLu2022,HuaLiDen2022}. In recent years, there is an increasing interest in applying deep learning techniques to learn from simulated physical systems. Several studies have explored this idea and proposed various methods for breaking down complex physical systems into simpler parts and analyzing their interactions to learn the underlying dynamics \cite{chang2017compositional,battaglia2016interaction}. Another possibility is to make use of Graph Neural Networks (GNN) on a complete graph to learn the dynamics by implicitly inferring the underlying interaction network \cite{Hoshen2017,Watters2017}. However, most of these previous works typically treat the interaction network recovery as the transient part of the learning process. In contrast, Neural Relational Inference (NRI) \cite{kipf2018_NRI} infers the interaction network explicitly, offering greater interpretability in the process. Hence, we employ neural relational inference to learn about the latent network edges and demonstrate it explicitly in several examples. We begin by discussing about NRI in the next section.

\begin{figure*}[!htbp]
    \centering
    \includegraphics[width=0.9\textwidth]{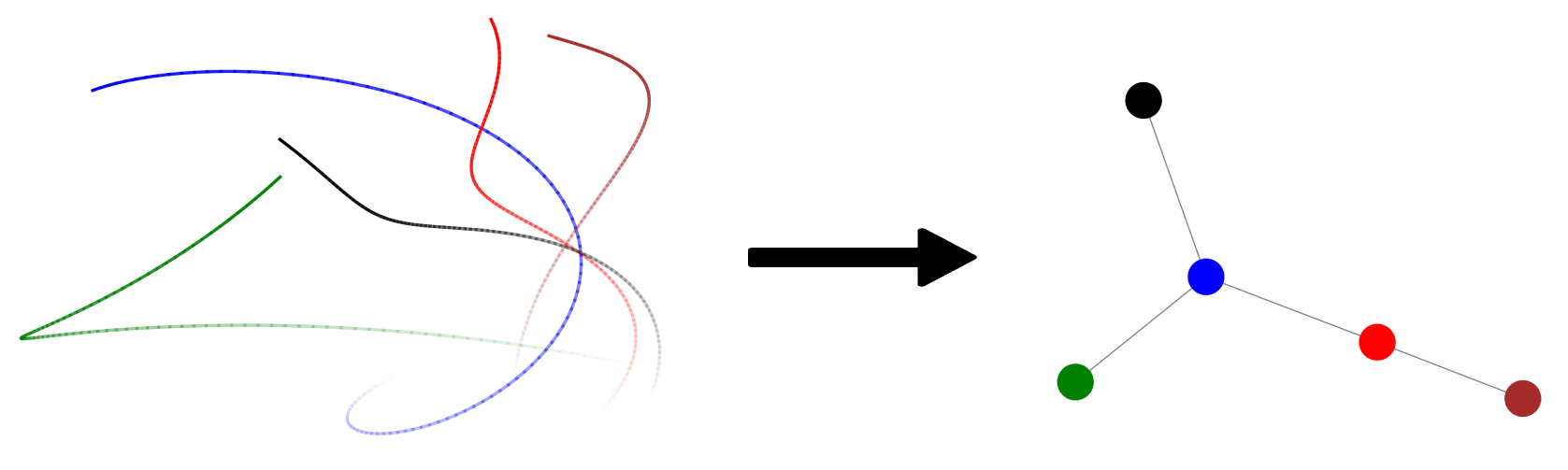}
    \caption{(left) Trajectories of 5 particles coupled via springs. (right) Corresponding interaction network is shown. This encodes the interaction between particles.}
    \label{fig:traj2net}
\end{figure*}

\section{Neural Relational Inference}
\label{nri_model}
Neural Relational Inference (NRI) is a self-supervised learning model that learns underlying interaction network and dynamics from observed data. NRI is a graph neural network-based variational autoencoder framework \cite{Kingma2014} in which the encoder learns the interaction network, and the decoder predicts the future trajectory of the an agent. In this work, we will adopt this learning model and we discuss its salient features in this section.
The NRI model uses GNN, which works by applying a message-passing algorithm on networks. As shown in Ref. \cite{gilmer2017neural}, one round of node-to-node message passing in GNN is defined as follows; given a graph $\mathscr{G} = (\mathscr{V}, \mathcal{E})$  with vertices $v \in \mathscr{V}$ and edges $e = (v, v') \in  \mathcal{E}$, for message passing from one node to another we have
\begin{equation} 
	v \rightarrow e :\textbf{m}_{(i,j)}^l = f_e^l ( \text{CONCAT} [ \textbf{m}_i^l,\textbf{m}_j^l,\textbf{x}_{(i,j)}]) ,
\label{msg_1}
\end{equation}
\begin{equation}
	e \rightarrow v : \textbf{ m}_{j}^{l+1} = f_v^l\left(\text{CONCAT}\left[ \sum_{i\in \mathscr{N}_j} \textbf{ m}_{(i,j)}^l, \textbf{x}_j \right]\right), 
\label{msg_2}
\end{equation}
where $ \textbf{m}_i^l $ is the $l$\textsuperscript{th} layer node embedding of $ v_{i} $ and $\textbf{m}_{(i,j)}^l$ is the edge embedding of $e_{(i,j)}$ at level $l$. Additionally, $\textbf{x}_i$ and $\textbf{x}_{(i,j)}$ refer to the node and edge features, respectively. Further, $\mathscr{N}_j$ denotes a set of node ids for nodes neighboring $j$-th node, while $\text{CONCAT}[.,.]$ signifies vector concatenation. The functions $f_v$ and $f_e$ denote neural networks specific to node and link, respectively. Equations \eqref{msg_1} and \eqref{msg_2} provide insight into how nodes and edges within the network communicate with one another.

The encoder employs two rounds of node to edge message passing to estimate the probability distribution of potential links $q_{\phi}(\textbf{z})$ using the trajectories of $\textbf{x} = (\textbf{x}_1^{1:T},...,\textbf{x}_N^{1:T})$ of $N$ agents taken as input. Then, the required probability distribution is given by
\begin{equation}
    q_{\phi}(\textbf{z}|\textbf{x}) = \text{softmax} \left( \textbf{m} ) \right) 
\end{equation}
where $\textbf{m} = f_{\text{ENC}}(\textbf{x})$ is a GNN action on a complete graph without self-loops.
As $q_{\phi}(\textbf{z}|\textbf{x})$ is a discrete distribution and hence not differentiable, the process of sampling links poses a challenge for backpropogation. To overcome this problem, NRI applies a continuous approximation of the discrete distribution to sample links making the process differentiable \cite{maddison2017concrete,jang2017categorical}. If $\textbf{z}_{ij}$ denotes the link type between the nodes $v_i, v_j$ and vector $\textbf{g} \in \mathbb{R}^K$ contains elements drawn from Gumbel(0, 1) distribution, then links are sampled in following manner:
\begin{equation}
    \textbf{z}_{ij} = \text{softmax} \left(\frac{\textbf{m}^2_{(i,j)} + \textbf{g}}{\tau}\right).
\end{equation}
Here $\textbf{m}^2_{(i,j)}$ denotes the edge embedding as obtained by the encoder after two rounds of node to edge message passing.
In this, the softmax temperature parameter $\tau$ controls the samples' ``smoothness". As $\tau \rightarrow 0$, the distribution becomes discrete.

The task of the decoder is to forecast the future evolution of dynamics using $p_{\theta}(\textbf{x}^{t+1}|\textbf{x}^t,...,\textbf{x}^1,\textbf{z})$ based on the initial state(s) of the system and the interaction network anticipated $\textbf{z}$ by the encoder. Decoder can be mathematically expressed as
\begin{equation}
    p_{\theta}(\textbf{x}|\textbf{z}) = \prod_{t=1}^{T}p_{\theta}(\textbf{x}^{t+1}|\textbf{x}^t,...,\textbf{x}^1,\textbf{z}).
\end{equation}
In many cases, calculating $p_{\theta}(\textbf{z}|\textbf{x})$ analytically using the Bayes rule (or, by other means)  is not feasible or computationally expensive due to difficulty in estimating $p_{\theta}(\textbf{x})$ (evidence) in high dimensional settings as is the case with the present model. Hence, we adopt an approximation
\begin{equation}
      p_{\theta}(\textbf{z}|\textbf{x}) \approx q_{\phi}(\textbf{z}|\textbf{x}).
\end{equation}
This approximation is effected by maximising Evidence Lower BOund (ELBO) given by \cite{Kingma2014}
\begin{equation}
    \mathcal{L} = \mathbb{E}_{q_\phi (\textbf{z}|\textbf{x})} \left[ \log p_{\theta}(\textbf{x}|\textbf{z}) \right] - \
    \text{KL} \left[ q_{\phi}(\textbf{z}|\textbf{x})||p_{\theta}(\textbf{z}) \right], 
\label{ELBO}
\end{equation}  
Where $p_{\theta}(\textbf{z})$ is the prior distribution over link types, which is assumed to be a uniform distribution, and $\text{KL}$ denotes the Kullback–Leibler divergence \cite{KL_Div} between two probability distributions. The additional details about ELBO are provided in the appendix \ref{appendix}.

\section{Dynamical Systems}
\label{dynsys}
In this section, we will apply this neural relational inference formalism for link inference and prediction of dynamics. We demonstrate the results with two different dynamical systems and their variants. We will call the first experiment ``Interacting Particles'' and the second ``Interacting Oscillators''. 

\subsection{Interacting Particles - Hooke's Law}
This system consists of $N$ particles (assumed to be point masses) trapped within a 2D box. A schematic of this system with $N=5$ is shown in Fig. \ref{fig:traj2net}. There is no external force on these particles except elastic collision with the walls of the box. We randomly couple (with probability 0.5) a pair of particles with a spring of spring constant $K = 1$. This results in the formation of an Erdos-Reyni (ER) type  network $G(n=N, p=0.5)$, where $p$ is the probability of a link characterised by $K=1$ (absence of an edge implies no spring). In this network, particles are the nodes and edges are the springs. The particles interact via Hooke's law $F_{ij} = -K(s_i - s_j)$, where $F_{ij}$ is the force applied by particle $v_j$ on particle $v_i$, whose position is denoted by $s_i$.

This approach can be extended to multiple link types as well. For instance, to simulate with three link types, we added an additional link type with $K = 0.5$ and all the three link types were sampled with equal probability. This can be viewed as generating an ER network $G\left(n=N, p_1 = \frac{1}{3},p_2 = \frac{1}{3}\right)$, where $p_1$ and $p_2$ are probabilities of coupling a node pair by a link of type $K=1$ and $K=0.5$, respectively. The initial position of particles are sampled from the normal distribution $\mathcal{N}(\mu = 0,\sigma = 0.5)$. For the initial velocity, a random vector with a magnitude of 0.5 is generated. Now, we can calculate the particle trajectories by solving Newton's equations of motion. The numerical solutions were obtained using leapfrog method with a step size of 0.001 \cite{Hairer_num_integ}. The input features to the model will be generated by concatenating the 2D position vector $(x,y)$ with the 2D velocity vector $(\Dot{x},\Dot{y})$. We then downsampled the features with a sampling frequency of 100 to get the final set of features (trajectories) for training and testing the model. In this experiment, as assumed in \cite{kipf2018_NRI}, though the particles interacted with each other, the boundary effects were almost negligible. But unlike Ref. \cite{kipf2018_NRI}, the conditions in which the boundary effects become as prominent as interparticle interactions are also explored here. This is done by changing the size of the 2D square box and also by changing its shape to circular.

\subsection{Interacting Oscillators - Kuramoto Model}
The Kuramoto oscillator model \cite{Kuramoto/BFb0013365} has been extensively used to study synchronisation in various physical settings ranging from biological (synchrony of neurons) to physical systems (Josephson arrays). This model describes a group of $N$ interacting oscillators where time evolution of the phase of i\textsuperscript{th} oscillator $\theta_i$ is given by:
\begin{equation}
\frac{d\theta_i}{dt} = \omega_i + \sum_{j \neq i} A_{ij} \sin(\theta_i - \theta_j) \label{K_ODE},
\end{equation}
where $\omega_i$ is the i\textsuperscript{th} oscillator's intrinsic frequency and $A_{ij}$ represents the coupling between oscillators indexed by $i$ and $j$.
The adjacency matrix $A_{ij}$ is also the interaction matrix and is taken to be an undirected and unweighted Erdos-Renyi network. Upon variation of coupling strength or the distribution of intrinsic frequencies, the Kuramoto model can display a variety of intricate dynamical behaviours including chaotic dynamics. In this work, $N=5$ and $N=10$ oscillators are employed and the initial values of intrinsic frequency and initial phase are sampled uniformly from intervals $[1,10)$ and $[0, 2\pi)$, respectively. Randomly chosen pairs $(i,j)$ of oscillators  are connected with coupling constant $A_{ij}$ such that $A_{ij}=1$ with probability $0.5$ and $A_{ij}=0$ with the same probability.

We have used the Dormand-Prince (DOPRI) method \cite{DORMAND198019} to solve Eq. \eqref{K_ODE} and obtain the time series of phases $\theta_i$. DOPRI is a Runge-Kutta method of order-5 with an embedded fourth order method for stepsize control. Further, $\theta_i$ is downsampled by a factor of 10 and is used to generate input features by concatenating $\frac{d\theta_i}{dt}$, $\sin\theta_i$ and  $\omega_i$. The input feature of the i\textsuperscript{th} oscillator is hence given by:

\begin{equation}
\textbf{x}_i = \text{CONCAT} \Bigl[ \frac{d\theta_i}{dt}, ~ \sin\theta_i, ~ \omega_i \Bigr], ~~~~~  i = 1, 2, \dots, N.
\end{equation}

However, note that $\omega_i$ is an intrinsic parameter of the Kuramoto model and is often inaccessible. In practice, the accessible information is limited to the time evolution of phases $\theta_i$. So unlike ref\cite{kipf2018_NRI}, where the actual frequency of each oscillator $\omega$ is given as one of the model inputs, we have given the estimated frequency of each oscillator $\omega_e$ instead of $\omega$. We estimated $\omega_e$  as the dc component of the Fourier transform of $\frac{d\theta_i}{dt}(t)$. Thus, in our case, $\omega$ need not be provided as an input to the problem. In the results, the performance of the model is compared for both scenarios, i.e., using actual $\omega$ and estimated $\omega_e$ for the input.

\section{Model performance}
The machine learning model described in Sec. \ref{nri_model} is employed to learn the connectivity matrix for the case of interacting particles and Kuramoto model described in Sec. \ref{dynsys}. In this section, we have evaluated the model on two different sets of tasks: interaction network recovery and trajectory forecasting. Performance on interaction network recovery tasks is quantified by accuracy in predicting the interaction type between the pair of agents, and performance in trajectory forecasting tasks is measured using the Mean Squared Error (MSE) in predicting the system's future state for $T$ timesteps.

The training data set comprises 50,000 simulations for all tasks mentioned in tables \ref{tab:1}-\ref{tab:4} and training is carried out for 500 epochs. To assess the model performance, testing is performed on a set of 10000 simulations for each task except for tasks 12 and 13 for which 2000 test simulations were used. The metrics presented in Table \ref{tab:1}-\ref{tab:4} were calculated based on the performance of the model on test data.

\subsection{Learning tasks : interaction network recovery}
\begin{figure}[!tbp]
    \centering
    \includegraphics[width=1\linewidth]{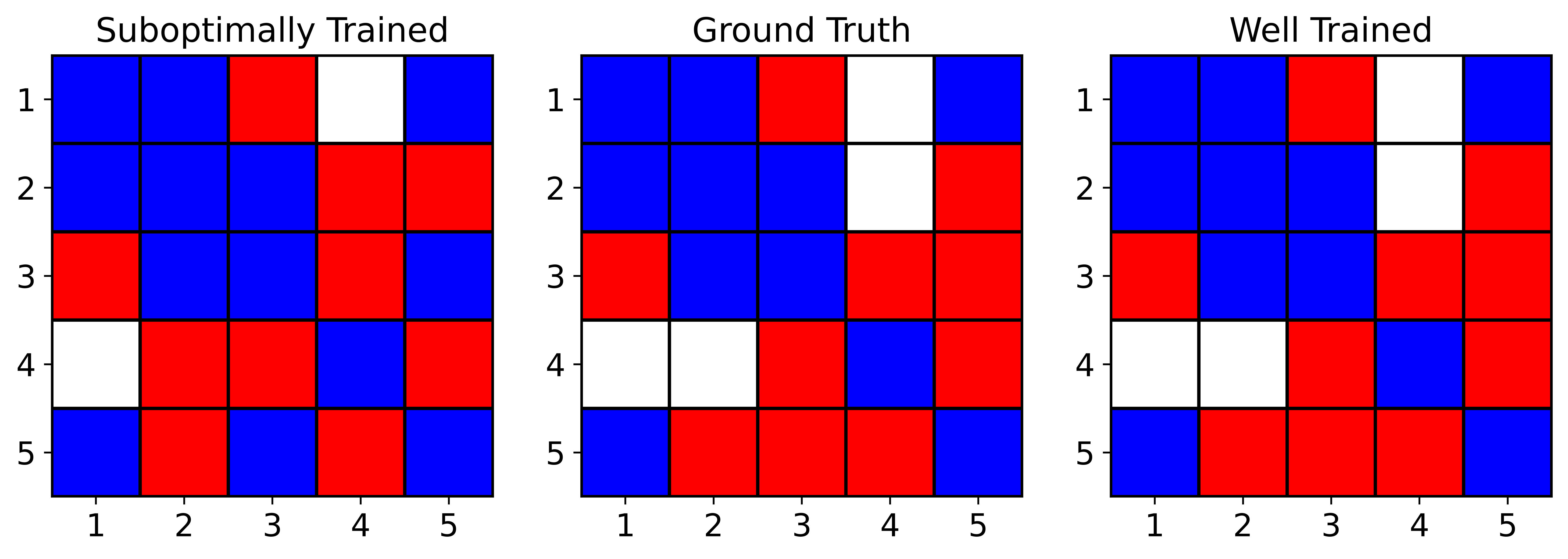}
    \caption{Performance comparison of suboptimally trained and well trained model in recovering the interaction network for task-3 in Table \ref{tab:1}. The numbers displayed on the axes of the adjacency matrices denote the oscillator index. (left) Adjacency matrix predicted by Suboptimally trained model. (right) Adjacency matrix predicted by well trained model. (middle) Actual (or, ground truth) adjacency matrix designed for the problem. Blue, red and white indicate that the particles are coupled by a link, respectively, of type $K=0$, $K=0.5$ and $K=1$. Predicted adjacency matrices are plotted for the majority class corresponding to the most probable link type anticipated by the encoder.}
    \label{fig:K_3 net}
\end{figure}

In this section, the results from interaction network recovery are presented. Tasks 1 and 2 (see Table \ref{tab:1}) consist of 5 and 10 interacting particles, respectively, where each pair of particles is either interacting (spring constant $K = 1$) or not interacting ($K = 0$) with equal probability. Here, we achieved an accuracy greater than 99\% with a small standard deviation, though the accuracy appears to decrease by about 1.5\% as the number of particles is doubled (task 2). Task 3 consists of 5 interacting particles where interaction is of 3 types -- can be of type $K =0, K = 0.5$ or $K = 1$ with equal probability. In this case too 99\% accuracy can be observed. Fig. \ref{fig:K_3 net} visualises the interaction network recovery performance of a model trained in task 3 on one test set simulation in two cases, one in which the model is trained suboptimally for just 50 epochs (suboptimally trained) and the other in which it is trained for complete 500 epochs (well trained). Fig. \ref{fig:K_3 net} reveals that even a suboptimally trained model can give reasonable accuracy, while a well trained model is almost entirely accurate.

\begin{figure}[!tbp]
    \centering
    \includegraphics[width=1\linewidth]{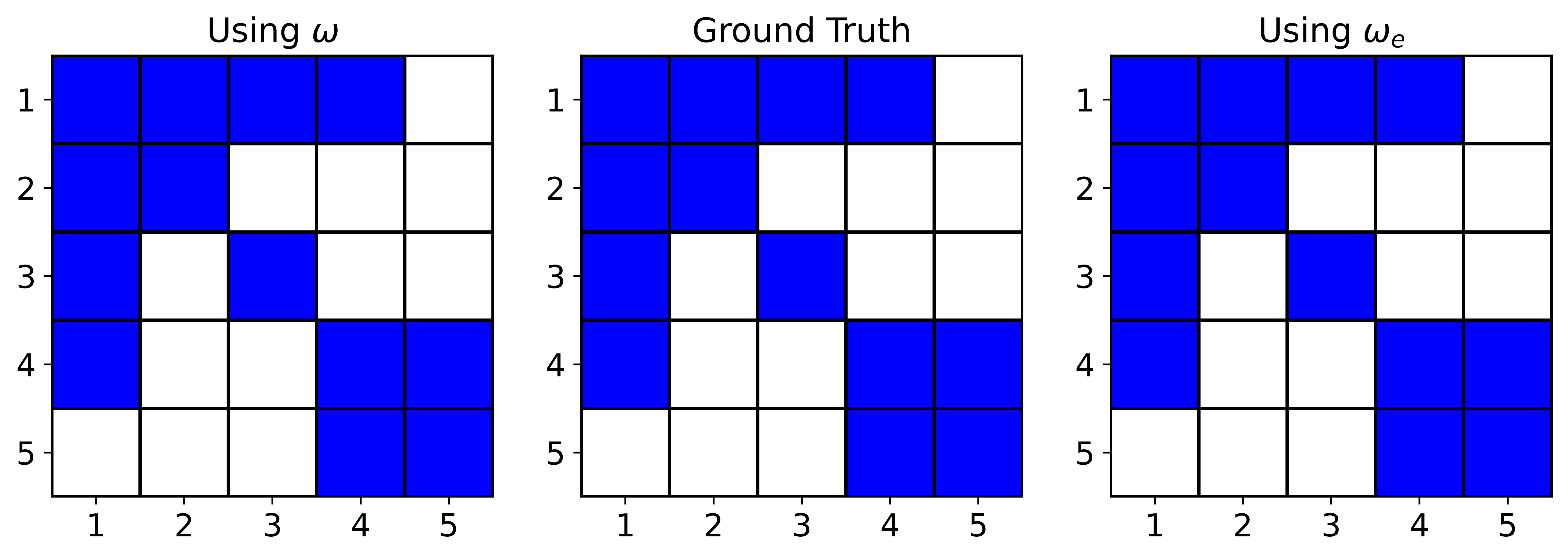}
\caption{Performance comparison of the model trained in task-4 (using the frequency $\omega$) as against the model trained in task-5 (using the estimate $\omega_e$) in recovering the interaction network. The numbers displayed on the axes of the adjacency matrices denote the oscillator index. White colur indicates that the oscillators are coupled, $A_{ij} = 1$, and blue indicates that they are not ($A_{ij} = 0$). (left) The predicted adjacency matrix when actual $\omega$ values are given as input to the model. (right) the predicted adjacency matrix when estimated frequency $\omega_e$ is given as input in lieu of actual values. (middle) actual (or, ground truth) adjacency matrix designed for the problem.}
    \label{fig:Kura_net}
\end{figure}

Tasks 4 to 7 consist of a system of 5 or 10 Kuramoto oscillators, which are either interacting or not interacting with equal probability. In tasks 4 and 6, actual values of the intrinsic frequency of the oscillators $\omega$ are given to the model as one of the inputs. In tasks 5 and 7, estimated intrinsic frequency $\omega_e$ is given in place of the actual value $\omega$. For a fixed number of oscillators, Table \ref{tab:1} shows that the accuracy obtained (if $\omega_e$ is used) is only mildly less than the case when $\omega$ is given as input. Fig. \ref{fig:Kura_net} shows interaction network recovery performance when $\omega$ is given, and $\omega_e$ is given for the case of 5 Kuramoto oscillators, confirms this conclusion. This implies that the machine learning model is robust since we need not know the values of characteristic frequencies $\omega$ to extract good predictive performance. However, as Table \ref{tab:1} reveals, an increase in the number of oscillators leads to substantially poor performance. While we expect this behaviour, a substantial decrease in accuracy is unexpected and requires to be investigated in detail.

Tasks 8 to 11 (Table \ref{tab:3}) are similar to task 1 (5 Particles, 2 link types) but with restrictive boundary conditions. In all other tasks pertaining to the dynamical system consisting of interacting particles (tasks 1-3, 12 and 13), the particles rarely interacted with the boundaries of the enclosing box due to the large box size and shorter trajectory length. However, in tasks 8 to 11, smaller box sizes and shapes were used for longer time duration of trajectories. This is done to assess the performance when boundary effects are as prominent as the interactions among agents. Typically, when boundary interaction begins to play a significant role as in the case of tasks 9 and 11 in Table \ref{tab:2}, the accuracy of interaction network recovery degrades substantially. As the boundary conditions become more restrictive, i.e, as the box size decreases, as seen in tasks 8-11 in Table \ref{tab:2}, the accuracy of the model in inferring interaction network also decreases. This is expected because now the dynamics is not only determined by interparticle interactions, but also by the elastic collision of particles with boundary walls adding to the system complexity. It is also evident that interaction network recovery is more difficult for a circular boundary than a square boundary.

In tasks 12 and 13 (Table \ref{tab:2}), the model trained in task 1 is evaluated in two extreme limits: (i) when the interaction network governing the dynamics is sparse ($G(n=5,p=0)$), and hence each particle is independent, and (ii) when the interaction network is dense ($G(n=5,p=1)$), and hence every particle is coupled to every other particle with $K = 1$. In both these cases, a mild decrease in accuracy is observed in comparison to task 1. As the model is trained on the simulation results in which the particles are either interacting or not interacting with probability 0.5, i.e., the interaction network is of the form $G(n=5,p=0.5)$, the model may have developed some bias towards the set of simulations used for training. However, this bias appears to be insignificant as the drop in accuracy is relatively low. This shows the generalisability of the model in wider settings. It is also worth noting that the drop in accuracy is more in the sparse regime (task 12) than the dense one (task 13). It indicates that the model is more inclined towards predicting the presence of a link in the interaction network rather than its absence.

\begin{table}[!tbp]
\centering
\begin{tabular}{|c|l|c|}
\hline \hline
\begin{tabular}[c]{@{}c@{}}Task\\      No.\end{tabular} & ~~~~~~~~~~~~Tasks                            & $\%$ Accuracy      \\ \hline
1                                                       & 5 Particles, 2 link types        & 99.868 $\pm$ 0.010 \\
2                                                       & 10 Particles, 2 link types       & 98.434 $\pm$ 0.034 \\
3                                                       & 5 Particles, 3 link types        & 99.171 $\pm$ 0.014 \\
\hline
4                                                       & 5 Oscillators ($\omega$ used)     & 96.340 $\pm$ 0.041 \\
5                                                       & 5 Oscillators ($\omega_e$ used)  & 95.270 $\pm$ 0.204 \\
6                                                       & 10 Oscillators ($\omega$ used)    & 74.648 $\pm$ 0.435 \\
7                                                       & 10 Oscillators ($\omega_e$ used) & 72.356 $\pm$ 2.678 \\ \hline \hline
\end{tabular}
\caption{Performance of the model on interaction network recovery tasks. Tasks 1-3 are for dynamics with Hooke's law type interactions. Tasks 4-7 are for Kuramoto oscillators for which either true frequency $\omega$ or estimated frequency $\omega_s$ is provided as input.}
\label{tab:1}
\end{table}

\begin{table}[b]
\centering
\begin{tabular}{|c|l|c|}
\hline \hline
\begin{tabular}[c]{@{}c@{}}Task\\      No.\end{tabular} & ~~~~~~~~~~~~~~~~~Tasks                           & $\%$ Accuracy      \\ \hline
8                                                      & Task 1 in square box of length = 4 & 95.964 $\pm$ 0.176 \\
9                                                      & Task 1 in square with length = 2 & 77.677 $\pm$ 1.402 \\
10                                                      & Task 1 in circular box of diameter = 4      & 94.430 $\pm$ 0.295 \\
11                                                      & Task 1 in circular box of diameter = 2      & 68.823 $\pm$ 1.673 \\  \hline
12                                                       & Task 1 with sparse interactions & 97.911 $\pm$ 1.951 \\
13                                                       & Task 1 with dense interactions & 99.785 $\pm$ 0.081 \\ \hline \hline
\end{tabular}
\caption{Further tests of interaction recovery tasks : Tasks 8 and 9 show the performance of the model trained in Task 1 (see Table \ref{tab:1}) when interaction network becomes sparse and dense. Tasks 10-13 are tests performed using task 1 with restricted boundary conditions.}
\label{tab:2}
\end{table}

\subsection{Learning tasks : Trajectory forecasting tasks}
\begin{figure*}[!htbp]
\centering
\begin{subfigure}{.5\textwidth}
  \centering
  \includegraphics[width=1\linewidth]{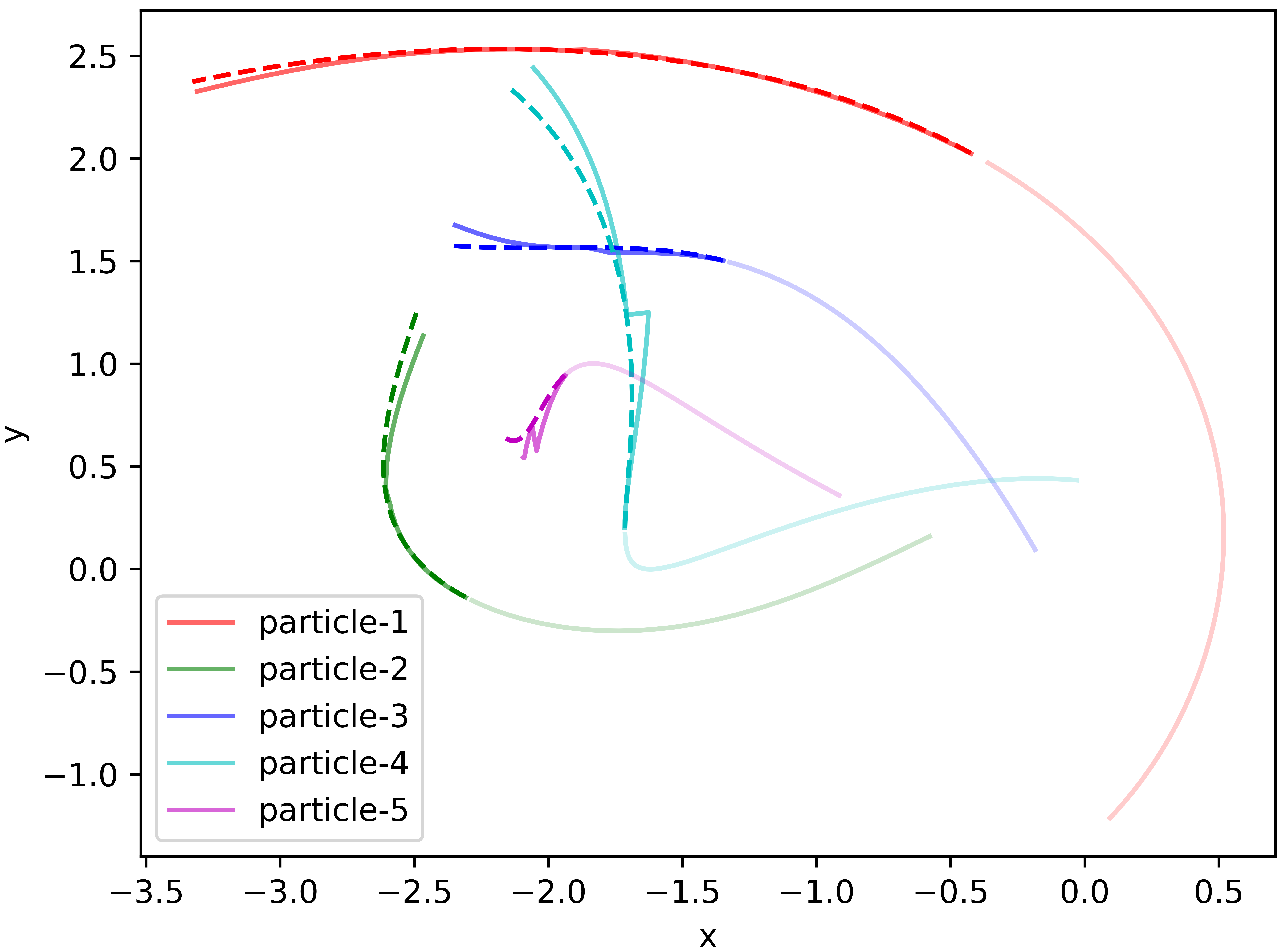}
  \caption{Suboptimally Trained}
  \label{fig:sub1}
\end{subfigure}%
\begin{subfigure}{.5\textwidth}
  \centering
  \includegraphics[width=1\linewidth]{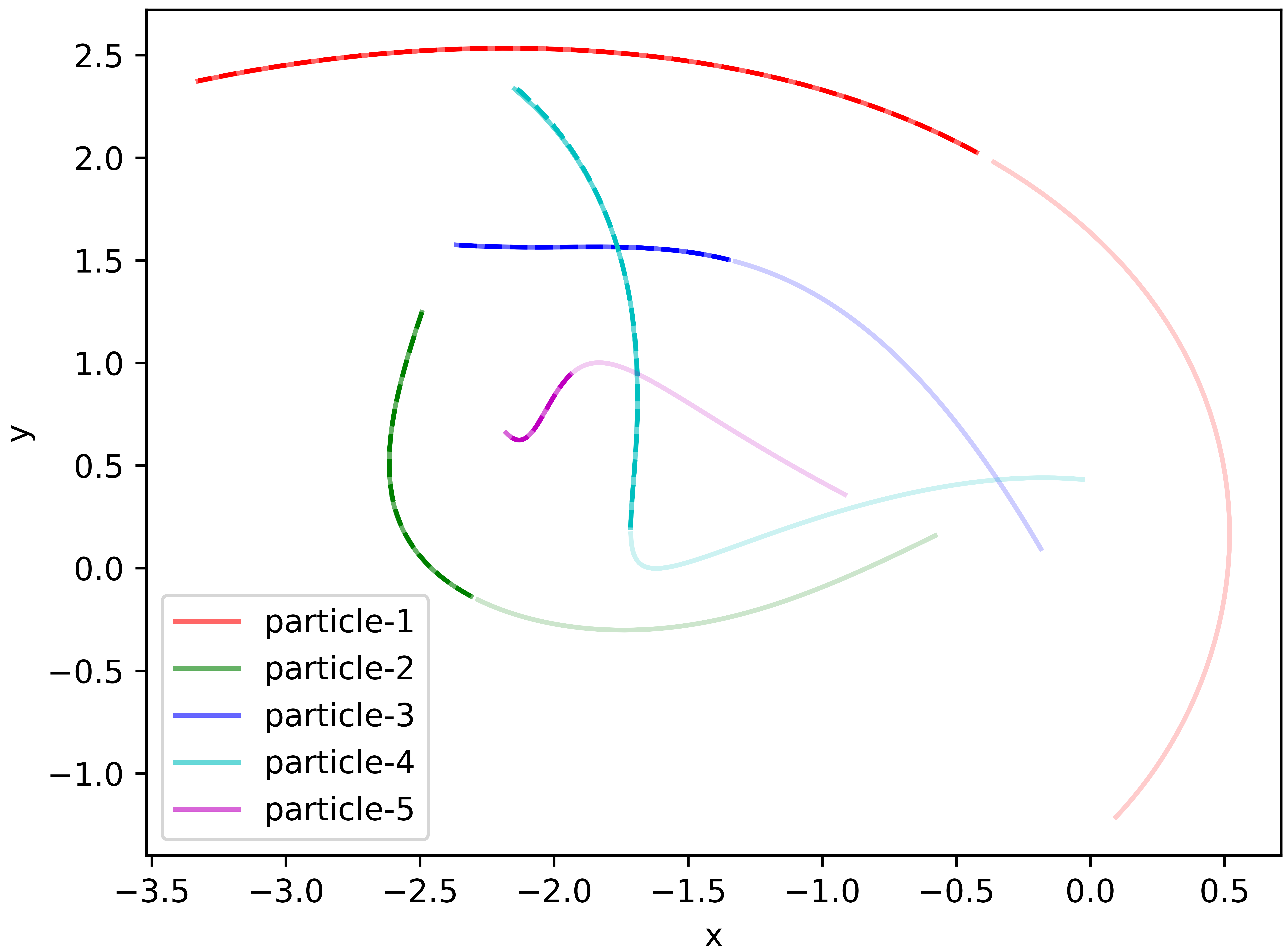}
  \caption{Well Trained}
  \label{fig:sub2}
\end{subfigure}
\caption{Performance comparison of (left) suboptimally trained model (trained for 10 epochs) and (right) well trained model (trained for 500 epochs) in forecasting the trajectory in task 3. The faint line of the trajectory corresponds to 49 timesteps given as input to the encoder, and the solid line denotes 40 future timesteps predicted by the decoder. The dashed line indicates the actual trajectory for 40 future timesteps.}
\label{fig:K3_traj}
\end{figure*}

Tables \ref{tab:3} and \ref{tab:4} shows the performance in predicting the future states of the dynamical systems for 10 and 20 timesteps. Model performance is evaluated using mean square error metric $\text{MSE} = \langle (\textbf{x}(t) - \textbf{x}_{\rm pred}(t))^2 \rangle$, where $\textbf{x}(t)$ represents the numerically exact trajectory and $\textbf{x}_{\rm pred}(t)$ is the trajectory predicted by the model. It is very clear from the tables that, for all the tasks, the MSE is relatively larger for $T = 20$ than $T = 10$. This is to be expected since predicting more steps into the future will accumulate more errors. Trajectory forecasting performance of suboptimally trained and well trained model in task 3 on a test simulation is shown in Fig \ref{fig:K3_traj}. This test simulation is the same one on which the performance of the interaction network recovery is visualized in Figure \ref{fig:K_3 net}. The Fig \ref{fig:K3_traj} indicates that a well trained model provides a better accuracy in the prediction of future dynamics. In general (as can be seen from Table \ref{tab:3}), MSE is less for the tasks involving the system of interacting particles than the system of oscillators. This may be because the system of oscillators shows much more complex and rich dynamics than the particles interacting via Hooke's law.

\begin{figure}[!htbp]
    \centering
    \includegraphics[width=1\linewidth]{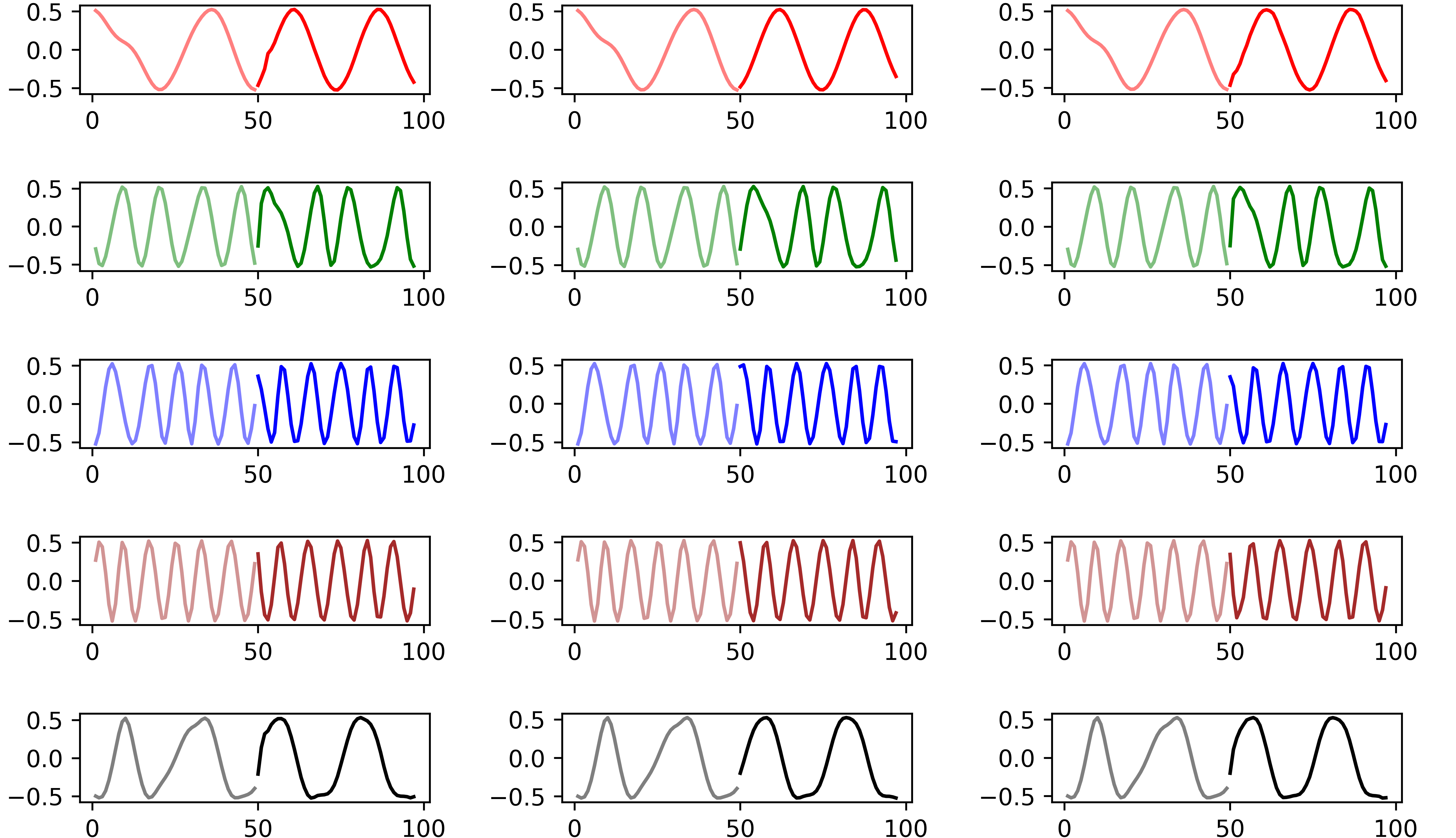}
    \subfloat[\label{actual} Prediction using $\omega$]{\hspace{.35\linewidth}}
    \subfloat[\label{ground} Ground Truth]{\hspace{.3\linewidth}}
    \subfloat[\label{actual} Prediction using $\omega_e$]{\hspace{.4\linewidth}}
    \caption{Performance comparison of the trajectories as a function of time trained in task 4 (using $\omega$, left panel) as against the model trained in task 5 (using estimated $\omega_e$, right panel) in forecasting the phase of oscillators. The middle panel is the trajectories from the numerical solution (the ground truth). Faint lines (trajectories) denote the first 49 timesteps given as model input, while solid trajectories denote the model predictions for the next 49 timesteps. Horizontal axes represent discrete time steps.}
    \label{fig:kura_traj}
\end{figure}

In the tasks 4-7 involving the system of oscillators in Table \ref{tab:3}, MSE is lower for the tasks where actual values of frequencies $\omega$  are given as inputs compared to the case when the estimated values of frequencies $\omega_e$ are given. This outcome can be anticipated since precise inputs facilitates better learning of the system dynamics. For the phases of Kuramoto oscillators, a detailed comparison of predictive performance on a test simulation (It is the same test simulation as used in Fig \ref{fig:Kura_net} to visualise the interaction network recovery performance) is shown for individual trajectories in Fig. \ref{fig:kura_traj}. It compares the predicted trajectories in both cases (when $\omega$ is given as one of the model inputs and when $\omega_e$ is given in place of $\omega$) with that of the actual (ground truth) trajectory. Since the errors involved are small, it is not clearly visible in the figure, though upon close examination, the results are consistent with that in Table \ref{tab:3}.

\begin{figure}[!htbp]
    \centering
    \includegraphics[width = 1\linewidth]{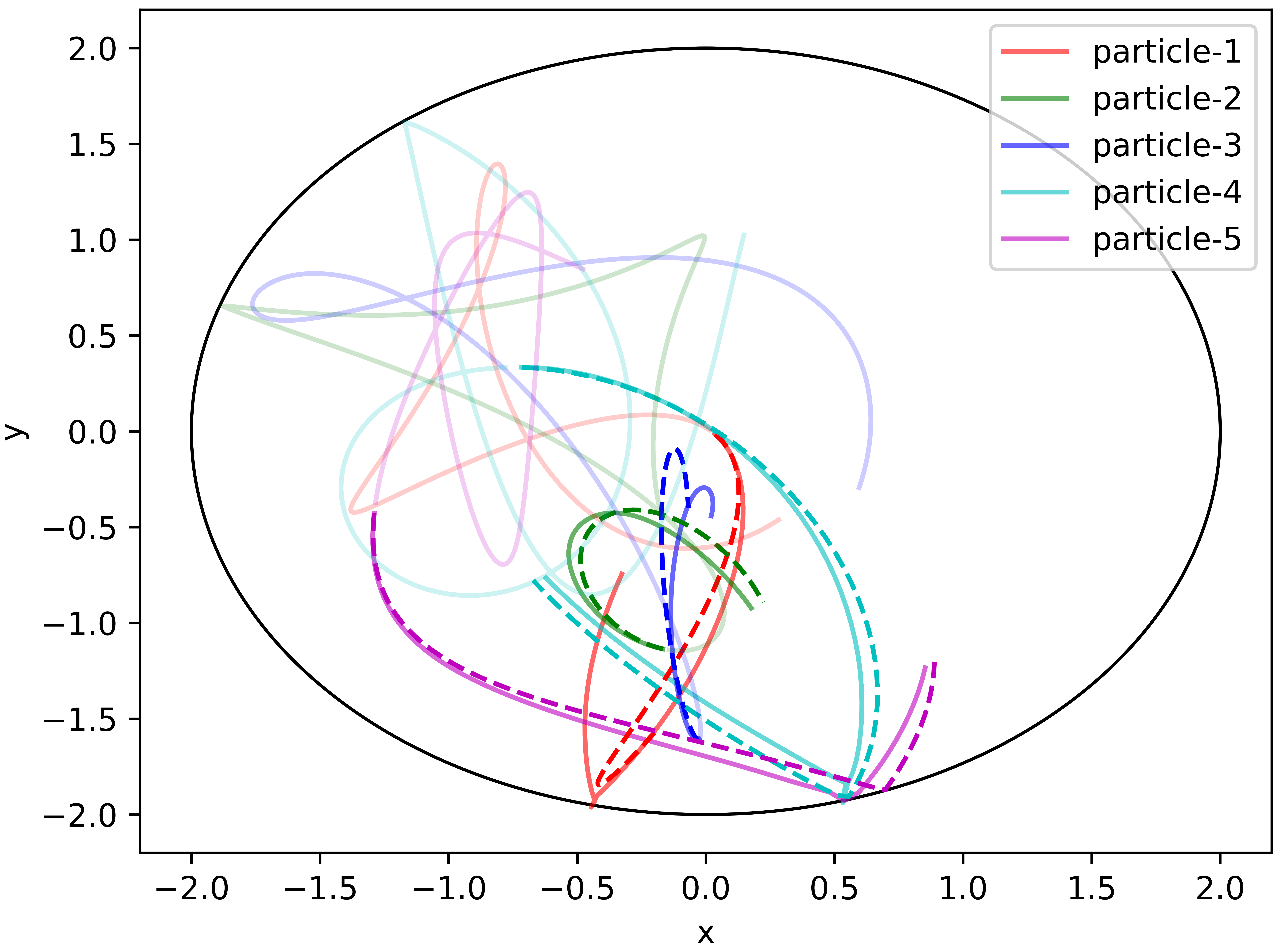}
    \caption{Performance of the model in forecasting trajectories under restrictive boundary conditions (particles are enclosed within a circle of radius 2, i.e. task 11). The faint part of the trajectory corresponds to the first 149 timesteps given as input to the encoder, and the solid part of the trajectory denotes the 60 future timesteps predicted by the decoder. The dashed part indicates the actual trajectory for 60 future timesteps, and the black circle denotes the circular boundary.}
    \label{fig:bcond2}
\end{figure}

As is evident from Table \ref{tab:4}, interactions with the boundary affects the trajectory predictions of the machine learning model. Similar to model performance in the interaction network recovery task, trajectory forecasting performance degrades with restrictive boundaries, i.e. MSE is seen to be increasing in tasks 8 to 11 as the size of the enclosure becomes smaller inducing more boundary encounters. A snapshot of trajectory predictions is shown for the case of task 11 in Fig. \ref{fig:bcond2}. The divergence of the solid trajectories (trajectories predicted by the model) from the dashed ones (actual trajectories) signifies the errors made by the model in forecasting future states of the dynamical system.

\begin{table}[!htbp]
\centering
\begin{tabular}{|c|l|c|c|}
\hline \hline
\multirow{2}{*}{\shortstack{Task \\ \\ No.}} & \multirow{2}{*}{~~~~~~~~Tasks} & \multicolumn{2}{c|}{MSE for $T$ time-steps} \\
  &                                  & $T = 10$ & $T = 20$ \\ \hline

1 & 5 particles, 2 link types        & $1.102 \times 10^{-6}$ & $4.580 \times 10^{-6}$ \\
2 & 10 particles, 2 link types       & $1.647 \times 10^{-5}$ & $5.138 \times 10^{-5}$ \\
3 & 5 particles, 3 link types        & $2.929 \times 10^{-6}$ & $1.057 \times 10^{-5}$ \\ \hline
4 & 5 oscillators (with $\omega$)     & $1.843 \times 10^{-3}$ & $4.258 \times 10^{-3}$ \\
5 & 5 oscillators (with $\omega_e$)  & $1.987 \times 10^{-3}$ & $5.238 \times 10^{-3}$ \\
6 & 10 oscillators (with $\omega$)    & $7.720 \times 10^{-3}$ & $2.217 \times 10^{-2}$ \\
7 & 10 oscillators (with $\omega_e$) & $7.996 \times 10^{-3}$ & $2.891 \times 10^{-2}$ \\ \hline \hline
\end{tabular}
\caption{Performance of the model on trajectory forecasting tasks for a variety of dynamical systems for $T=10$ and $T=20$ time steps. Performance measured using mean squared error metric.}
\label{tab:3}
\end{table}

\begin{table}[!htbp]
\begin{tabular}{|c|c|c|c|}
\hline \hline
\multirow{2}{*}{\shortstack{Task \\ \\ No.}}  & \multirow{2}{*}{Tasks} & \multicolumn{2}{c|}{MSE for $T$ time-steps} \\
  &                                  & $T = 10$ & $T = 20$ \\ \hline
8 & Square with length = 4 &   $1.697 \times 10^{-4}$       &  $4.886 \times 10^{-4}$       \\
9 & Square with length = 2 &$ 1.339 \times 10^{-3} $      & $4.610 \times 10^{-3}$       \\
10 & Circle with diameter = 4      & $3.036 \times 10^{-4}$       & $8.047 \times 10^{-4}$       \\
11 & Circle with diameter = 2      & $1.298 \times 10^{-3}$       & $4.234 \times 10^{-3} $      \\ \hline \hline
\end{tabular}
\caption{Performance of the model in forecasting future trajectories in case of 5 Particles, 2 link
types subjected to the restrictive boundary conditions. Performance measured using mean squared error metric. }
\label{tab:4}
\end{table}

\section{Summary and Outlook}
In this work, Neural Relational Inference, a graph neural network-based self-supervised learning model has been employed to infer interaction network and trajectories for two distinct many-body coupled dynamical systems. One of these is a linear system, namely, particles coupled through Hooke's law type force. The other system is a nonlinear system, namely, a collection of coupled Kuramoto oscillators. It is demonstrated that the neural relational inference model has performed well at predicting the interaction matrix that encodes information about the connectivity structure between the particles. Further, the same machine learning model also provides a good prediction for trajectories of the dynamical systems. Remarkably, the model does not make any assumptions regarding the dynamical system, especially about the type of interaction, and this indicates the generality of this approach. 

In the case of the Kuramoto oscillator system, it is also shown that instead of giving the actual value of intrinsic frequency (often an inaccessible system parameter) as input to the model, a frequency value estimated as the DC component of the Fourier transform of the trajectory also has good predictive power. This approach makes the model agnostic to system parameters. Further, the applicability of the NRI model has been extended by successfully applying it to dynamical systems where boundary effects are as prominent as interactions among agents. It is shown that as boundary effects become prominent, the model performance degrades on both fronts, namely interaction network recovery and trajectory forecasting.

There is significant room for improvements as Interaction Network Inference is still in its early stage of development. The NRI model is less effective when the number of particles are large and hence the interaction networks are also large. This is primarily because the resource and time complexity of the model is $\mathcal{O}(|\mathcal{V}|^2)$ where $\mathcal{V}$ denotes a set of nodes in the interaction network. The NRI model starts with a fully connected network and discovers the sparsity gradually. This might be unnecessary since most real-world networks are sparse. Improvements in these directions provide fruitful avenues for future research. Another possible direction in the NRI model is to extend it to the cases where the interaction network itself is time-dependent. We have only considered two deterministic dynamical systems, another challenge would be to extend this framework to predict interaction network and trajectories of stochastic dynamical systems.\newline

\section*{Acknowledgements}
We acknowledge the High-Performance Computing (HPC) facility, Param Brahma, at the Indian Institute of Science Education and Research (IISER) Pune, for providing the computational resources required for this research.

%

\clearpage
\onecolumngrid
\appendix
\appendix

\section{Evidence Lower Bound (ELBO)}
\label{appendix}

By  Bayes rule, posterior distribution 
  \begin{equation}
      p_{\theta}(\textbf{z}|\textbf{x}) = \frac{p_{\theta}(\textbf{x}|\textbf{z}) p_{\theta}(\textbf{z})}{p_{\theta}(\textbf{x})}
  \end{equation}

 In many cases (as is the case with the NRI model), calculating $p_{\theta}(\textbf{z}|\textbf{x})$ analytically (e.g. using the Bayes rule)  is not possible or computationally expensive (infeasible) due to difficulty in estimating $p_{\theta}(\textbf{x})$ (evidence) in high dimensional settings. So we approximate $p_{\theta}(\textbf{z}|\textbf{x})$ by $q_{\phi}(\textbf{z}|\textbf{x})$ i.e.
\begin{equation}
     q_{\phi}(\textbf{z}|\textbf{x}) \approx p_{\theta}(\textbf{z}|\textbf{x}) 
\end{equation}
\textbf{Claim}: Minimising $KL(q_{\phi}||p_{\theta})$ is equivalent to maximising ELBO.
    \begin{proof}
  
  \begin{align*}
  KL(q_{\phi}||p_{\theta}) &= \mathbb{E}_{q_{\phi}}\left[\text{log} \frac{q_{\phi}(\textbf{z}|\textbf{x})}{p_{\theta}(\textbf{z}|\textbf{x})}\right]\\
  &= \mathbb{E}_{q_{\phi}}[\text{log} q_{\phi}(\textbf{z}|\textbf{x})]-\mathbb{E}_{q_{\phi}}\left[\text{log} \frac{p_{\theta}(\textbf{z},\textbf{x})}{p_{\theta}(\textbf{x})}\right] \\
  &= \mathbb{E}_{q_{\phi}}[\text{log} q_{\phi}(\textbf{z}|\textbf{x})]-\mathbb{E}_{q_{\phi}}[\text{log} p_{\theta}(\textbf{z},\textbf{x})] + \text{log} p_{\theta}(\textbf{x}) \int q_{\phi}(\textbf{z}|\textbf{x}) dz \\
  &= \mathbb{E}_{q_{\phi}}[\text{log} q_{\phi}(\textbf{z}|\textbf{x})]-\mathbb{E}_{q_{\phi}}[\text{log} p_{\theta}(\textbf{z},\textbf{x})] + \text{log} p_{\theta}(\textbf{x})  \\
  \implies \underbrace{\text{log} p_{\theta}(\textbf{x})}_\text{log evidence}  &= \underbrace{-\mathbb{E}_{q_{\phi}}[\text{log} q_{\phi}(\textbf{z}|\textbf{x})] + \mathbb{E}_{q_{\phi}}[\text{log} p_{\theta}(\textbf{z},\textbf{x})]}_\text{term 1} + \underbrace{KL(q_{\phi}||p_{\theta})}_\text{term 2 $\ge$ 0} \\
  \implies \underbrace{\text{log} p_{\theta}(\textbf{x})}_\text{log evidence}  &\ge \underbrace{-\mathbb{E}_{q_{\phi}}[\text{log} q_{\phi}(\textbf{z}|\textbf{x})] + \mathbb{E}_{q_{\phi}}[\text{log} p_{\theta}(\textbf{z},\textbf{x})]}_\text{Evidence Lower Bound (ELBO)}  
  \end{align*}   
\end{proof}
In the penultimate step, it is evident that minimising term 2 is equivalent to maximising term 1 to maintain equality.
In the last step, we arrived at an inequality by virtue of the fact that KL divergence is positive semidefinite (can be proved using Jensen's inequality and by noting that log is a concave function). The last step also justifies calling RHS of inequality as Evidence Lower Bound (ELBO).\par

We can easily transform ELBO into the form used in the main body of the paper.
\begin{equation}
\begin{aligned}
    \text{ELBO} &= -\mathbb{E}_{q_{\phi}}[\text{log} q_{\phi}(\textbf{z}|\textbf{x})] + \mathbb{E}_{q_{\phi}}[\text{log} p_{\theta}(\textbf{z},\textbf{x})]\\
    &= -\mathbb{E}_{q_{\phi}}[\text{log} q_{\phi}(\textbf{z}|\textbf{x})] + \mathbb{E}_{q_{\phi}}[\text{log} p_{\theta}(\textbf{x}|\textbf{z})] + \mathbb{E}_{q_{\phi}}[\text{log} p_{\theta}(\textbf{z})]\\
    &= \mathbb{E}_{q_{\phi}}[\text{log} q_{\phi}(\textbf{x}|\textbf{z})] - \mathbb{E}_{q_{\phi}}  \left[\text{log} \frac{q_{\phi}(\textbf{z}|\textbf{x})}{p_{\theta}(\textbf{z})}\right]
\end{aligned}
\end{equation}
\begin{equation}
    \text{ELBO} = \mathbb{E}_{q_\phi (\textbf{z}|\textbf{x})}[log p_{\theta}(\textbf{x}|\textbf{z})] - KL[q_{\phi}(\textbf{z}|\textbf{x})||p_{\theta}(\textbf{z})]
\end{equation}

Maximising ELBO means maximising $\mathbb{E}_{q_\phi (\textbf{z}|\textbf{x})}[log p_{\theta}(\textbf{x}|\textbf{z})]$ and minimising $KL[q_{\phi}(\textbf{z}|\textbf{x})||p_{\theta}(\textbf{z})]$. Maximising $E_{q_\phi (\textbf{z}|\textbf{x})}[log p_{\theta}(\textbf{x}|\textbf{z})]$ is mathematically equivalent to minimising reconstruction error.  The term $KL[q_{\phi}(\textbf{z}|\textbf{x})||p_{\theta}(\textbf{z})]$ acts as a regularisation term.

\end{document}